\documentstyle{article}                   
 \pdfoutput=1 \pdfcompresslevel=9 \pdfadjustspacing=1
\def\ifudf#1{\expandafter\ifx\csname #1\endcsname\relax}
\newif\ifpdf \ifudf{pdfoutput}\pdffalse\else\pdftrue\fi
\everydisplay={\textstyle}
\font\tbb=bbmsl10 \font\sbb=bbmsl10 scaled 700
 \font\fbb=bbmsl10 scaled 500
\newfam\bbfam \textfont\bbfam=\tbb
 \scriptfont\bbfam=\sbb \scriptscriptfont\bbfam=\fbb
\def\bb{\fam\bbfam}%
\font\tfk=eufm10 \font\sfk=eufm7 \font\ffk=eufm5
\newfam\fkfam \textfont\fkfam=\tfk
 \scriptfont\fkfam=\sfk \scriptscriptfont\fkfam=\ffk
\def\fn[#1]{\font\TmpFnt=#1\relax\TmpFnt\ignorespaces}
\def\em{\expandafter\ifx\the\font\tensl\rm\else\sl\fi}
\def\dt{\number\day.\number\month.\number\year/\the\time}
\def\\{\hfill\break}
\newcount\secNo \secNo=0
\def\section #1\par{\goodbreak\vskip 3ex\noindent
 \global\advance\secNo by 1 \global\eqnNo=0
 {\fn[cmbx10 scaled 1200]#1}\vglue 1ex}
\def\href#1<#2>{\leavevmode
 \ifpdf\pdfstartlink attr {/Border [0 0 0 ]} goto name {#1}\fi
 {#2}\ifpdf\pdfendlink\fi}
\def\label@#1:#2@{\ifudf{#1}
 \expandafter\xdef\csname#1\endcsname{#2}\else
 \errmessage{label #1 already in use!}\fi}
\label@fig.ribs:Fig 1@
\label@tcr:2.1@
\label@elift:2.2@
\label@tcircle:2.3@
\label@ecircle:2.4@
\label@fig.ribd:Fig 2@
\label@plift:2.5@
\label@csphere:2.6@
\label@nisomc:2.7@
\label@thm.pnf1:Lemma 2.1@
\label@psphere:2.8@
\label@thm.pnf2:Lemma 2.2@
\label@thm.rch:Cor 2.3@
\label@thm.sch:Cor 2.4@
\label@fig.chls:Fig 3@
\label@drec:2.9@
\label@thm.sdr:Cor 2.5@
\label@fig.redd:Fig 4@
\label@thm.ddr:Thm 2.6@
\label@thm.rc1:Cor \& Def 2.7@
\label@fig.intd:Fig 5@
\label@demoulin:2.10@
\label@thm.ri2:Thm 2.8@
\label@thm.ri3:Cor 2.9@
\label@def.ribs:Def 3.1@
\label@thm.ribt:Lemma 3.2@
\label@derb:3.1@
\label@def.ribd:Def 3.3@
\label@nisomm:3.2@
\label@thm.pnfm:Lemma 3.4@
\label@thm.rchm:Cor 3.5@
\label@thm.sdrm:Cor 3.6@
\label@drem:3.3@
\label@thm.ddrm:Cor 3.7@
\label@permij:3.4@
\label@bquad:3.5@
\label@thm.src:Thm \& Def 3.8@
\label@thm.drc:Thm \& Def 3.9@
\label@imdg:10@
\def\@#1:#2@{\ifpdf\pdfdest name {#1} xyz\fi {#2}}
\def\:#1:{\href#1<\ifudf{#1}??\else\csname#1\endcsname\fi>}
\newcount\refNo \refNo=0
\def\refitem#1 {\global\advance\refNo by 1
 \item{\@#1:\number\refNo@.}}
\let\pleqno\eqno \newcount\eqnNo \eqnNo=0
\def\eqno#1$${\global\advance\eqnNo by 1
 \pleqno{\rm(\@#1:\number\secNo.\number\eqnNo@)}$$}
\def\Z{{\bb Z}} \def\R{{\bb R}} \def\H{{\bb H}} \def\P{{\bb P}}
\def\fdeb#1{#1^{(1)}} \def\wlnb#1{#1^{(1)\perp}/#1}
 \def\Im{\mathop{\rm Im}}
\newtoks\title \newtoks\stitle \newtoks\author
\newtoks\status \newtoks\funding
\title={Ribaucour coordinates}
\stitle=\title
\author={F Burstall, U Hertrich-Jeromin, M Lara Miro}
\status={F Burstall et al}
\funding={Furthermore,
 we gratefully acknowledge financial support
 of the second author by the 
  Bath Institute for Mathematical Innovation;
 and
 of the third author by the
  TU Wien Doctoral Training Centre ``Computational Design''.
 }
\ifpdf\fi
\def\item#1{\par\leavevmode\hangindent=\parindent\hangafter=1%
 \llap{#1\enspace}\ignorespaces}
\begin{document}                          
\centerline{{\fn[cmbx10 scaled 1440]\the\title}}\vglue .2ex
\centerline{{\fn[cmr7]\the\author}}\vglue 3em plus 3ex
\centerline{\vtop{\hsize=.8\hsize{\bf Abstract.}\enspace
 We discuss results for the Ribaucour transformation of curves
 or of higher dimensional smooth and discrete submanifolds.
 In particular, a result for the reduction of the ambient
 dimension of a submanifold is proved and
 the notion of Ribaucour coordinates is derived using a
 Bianchi permutability result.
 Further, we discuss smoothing of semi-discrete curvature
 line nets and an interpolation by Ribaucour transformations.
}}\vglue 2em
\centerline{\vtop{\hsize=.8\hsize{\bf MSC 2010.}\enspace
 {\it 53C42\/}, {\it 53A10\/}, 53A30, 37K25, 37K35
}}\vglue 1em
\centerline{\vtop{\hsize=.8\hsize{\bf Keywords.}\enspace
 channel surface; canal surface; semi-discrete surface;
 curvature line net; discrete principal net; circular net;
 Ribaucour transformation; Ribaucour coordinates.
}}\vglue 3em

\section 1. Introduction

This paper touches upon several ideas and results concerning
the Ribaucour transformation of submanifolds and,
in particular, of curves:
the $1$-dimensional case not only helps to illustrate
the main ideas of our investigations, but it is also of
interest for the construction of $2$-dimensional discrete
and semi-discrete principal ``circular'' nets and for their
``smoothing'',
cf [\:bjmr16:].
On the other hand, our discussion of higher dimensional
submanifolds not only provides a generalization of
some of the results for curves, but it also sheds light
on the reasons and structure behind the results for curves.
Thus we first present several results for the Ribaucour
transformation of curves and then discuss whether and
how these generalize to higher dimensions:
a reduction of the ambient dimension by means of Ribaucour
transformation, which leads to Ribaucour coordinates;
smoothing of a sequence of Ribaucour transforms;
and
the possibility of an interpolation between curves or
submanifolds by means of a sequence of Ribaucour
transformations.
Each of these problems is addressed in both
the smooth and discrete settings.

Ribaucour coordinates were used by the classical authors
to investigate the geometry of surfaces with particular
properties of their curvature lines:
given a surface and a suitable plane,
for example, a tangent plane of the surface,
one may locally construct a regular map of the surface
to the plane by means of touching $2$-spheres.
These ``Ribaucour coordinates'' then map the curvature line
net of the surface to an orthogonal net in the plane.
For example,
it is advantageous to employ this type of coordinates to
investigate and construct surfaces of constant mean curvature
$H=1$ in hyperbolic space,
cf [\:liro02:] or [\:imdg:,\S5.5.27].
These classical Ribaucour coordinates were generalized in
more recent work to hypersurfaces,
see [\:cote04:, Cor 2.10],
and to submanifolds with flat normal bundle,
see [\:dft07:, Thm 1].

One principal aim of this paper is to provide a more direct
geometric approach to these generalized Ribaucour coordinates,
much in the spirit of the classical authors.
This approach relies on a geometric method to reduce the
ambient dimension of a curve or submanifold by means of
a Ribaucour transformation into a hypersphere,
see \:thm.sdr:, \:thm.ddr: for (smooth and discrete) curves
and \:thm.sdrm:, \:thm.ddrm: for submanifolds.
Once the dimension reduction is established,
the higher codimension version of Ribaucour coordinates
relies on a Bianchi permutability result,
cf [\:dft07:, Thm 2] and [\:imdg:, \S8.5.8].
In our setting, this permutability result is readily verified
by elementary means,
leading to Ribaucour coordinates for any smooth submanifold
with flat normal bundle resp any discrete circular net
of higher codimension,
see \:thm.src: and \:thm.drc:.

Another observation that relies on rather similar ideas
leads to a ``smoothing'' of a semi-discrete curvature
line net,
since such a net can be thought of a a sequence of Ribaucour
transforms of one curve.
If two $m$-dimensional submanifolds envelop an $m$-sphere
congruence then this sphere congruence provides a metric
and connection preserving isomorphism field between normal
bundles,
see \:thm.pnf1: and \:thm.pnfm:,
cf [\:cote04:, Cor 2.9].
As a consequence, any pair of curves that envelop
a $1$-parameter family of circles,
that is, form a Ribaucour pair,
are two curvature lines of a channel surface,
see \:thm.sch: or \:fig.chls::
this result can be thought of as the semi-discrete
version of a ``smoothing'' result for discrete circular
nets [\:bhv12:] and [\:bkpww11:],
where a discrete circular net is ``smoothed'' by fitting
Dupin cyclide patches into its facets.
The higher dimensional version \:thm.rchm: emphasizes that
the construction only depends on an enveloped $m$-sphere
congruence, not the fact that the two submanifolds form
a Ribaucour pair:
if they do, more structure can be obtained.

Finally, we prove that any two (smooth or discrete) curves
can be transformed into each other by a sequence of Ribaucour
transforms, see \:thm.ri3:.
In contrast to Cauchy problems for discrete nets,
cf [\:boje01:, Sect 4] or [\:imdg:, \S8.4.9],
this is a mixed boundary and initial value problem
and allows to construct discrete or semi-discrete principal
nets with given boundary data.
This result hinges again on the aforementioned dimension
reduction,
which allows us to reduce the problem to planar curves.
For planar curves, the result is then proved
by a simple geometric construction in the discrete case,
and as a trivial consequence of [\:buje06:, Thm 3.4]
 in the smooth case, see \:thm.ri2:.
However, as the flatness of a certain vector bundle
is required
---
which is trivial in the $1$-dimensional case,
but not in higher dimensions
---
our interpolation result does not generalize to higher
dimensions.

Though we start by formulating our results for curves in
Euclidean space in order to make the text more accessible,
we quickly resort to sphere geometric methods:
as the Ribaucour transformation is a Lie geometric notion,
many arguments are more efficient and transparent in the
M\"obius or Lie geometric settings.
Since these techniques are discussed in great detail in other
works, we will only give few hints or details where required,
and refer the reader to, for example, [\:imdg:]
as well as to our previous papers [\:buje06:] or [\:bjmr16:].

{\it Acknowledgements.\/}
This work would not have been possible without the valuable
and enjoyable discussions with
C M\"uller,
M Pember,
F Rist
and
G Szewieczek
about the subject.

\the\funding

\section 2. Space curves

To set the scene we discuss smooth and discrete curves
in $3$-space:
this case already displays most of the key ideas of this
study, without the additional complexity caused by the
partial differential equations occurring in the case
of higher dimensional submanifolds.

 \pdfximage width 200pt {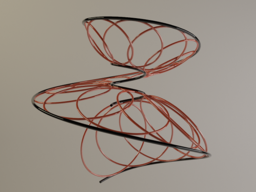}
  \xdef\ribsA{\the\pdflastximage}
 \pdfximage width 200pt {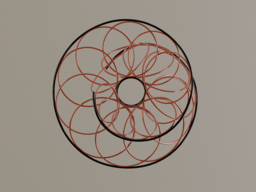}
  \xdef\ribsB{\the\pdflastximage}
 \vglue 0pt\hbox to \hsize{\hfil
  \vtop{\pdfrefximage\ribsA}%
  \hfil
  \vtop{\pdfrefximage\ribsB}%
 \hfil}\vglue 6pt\centerline{%
  {\bf\@fig.ribs:Fig 1@.}
   A Ribaucour pair of curves}\par

Thus let $x,\hat x:I\to\R^3$ denote a {\em Ribaucour pair\/}
of curves,
that is,
$x$ and $\hat x$ share a common tangent circle $c(u)$
at corresponding points $x(u)$ and $\hat x(u)$,
cf [\:bjmr16:, Def 2.1] and \:fig.ribs:.
We say that the two curves are mutual
{\em Ribaucour transforms\/} of each other.
For regularity we assume $x(u)\neq\hat x(u)$ for all
$u\in I$.
Analytically, this relation between the two curves can be
encoded by the reality of their {\em tangent cross ratio\/},
$$
  cr := x'(\hat x-x)^{-1}\hat x'(\hat x-x)^{-1}:I\to\R,
\eqno tcr$$
where products and inverses of vectors $v\in\R^3\cong\Im\H$
are computed in terms of Clifford, or quaternionic,
multiplication,
see [\:bjmr16:, Lemma 2.2].
This (quadratic) relation between the two curves may be
linearized by means of standard M\"obius geometric methods,
see [\:bjmr16:, Lemma 2.3]:
replacing $\R^3$ by the conformal $3$-sphere,
thought of as the projective light cone in $\R^{4,1}$,
$$
  S^3 \cong \P({\cal L}^4)
   \enspace{\rm with}\enspace
  {\cal L}^4 = \{y\in\R^{4,1}\,|\,(y,y)=0\},
$$
the circularity condition becomes linear in terms
of homogeneous coordinates $x$ and $\hat x$.
In detail, consider the orthogonal decomposition
$$
  \R^{4,1} = \R^{1,1}\oplus_\perp\R^3,
   \enspace{\rm where}\enspace
  \R^{1,1} = \langle o,q\rangle
$$
is spanned by isotropic vectors $o,q\in\R^{4,1}$ with
inner product $(o,q)=-1$,
and the {\em Euclidean lift\/} of a curve in $\R^3$ into
the light cone obtained from the isometric embedding
$$
  \R^3\ni x\mapsto\xi := o+x+{1\over 2}(x,x)\,q
   \in {\cal Q}^3 := \{y\in{\cal L}^4\,|\,(y,q)=-1\}
   \subset{\cal L}^4.
\eqno elift$$
The circularity condition of the Ribaucour transformation
for two curves now reduces to linear dependence:
$
  u\mapsto\dim\langle\xi,\xi',\hat\xi,\hat\xi'\rangle(u)=3.
$
Clearly, this condition is independent of the scaling of
{\em lifts\/}
$\xi\in\Gamma(\langle o+x+{1\over 2}(x,x)\,q\rangle)$ of $x$
and $\hat\xi$ of $\hat x$,
in fact, the map
$$
  u\mapsto c(u) := \langle\xi,\xi',\hat\xi,\hat\xi'\rangle(u)
\eqno tcircle$$
is independent of lifts and, in the case of a Ribaucour pair,
encodes the enveloped circle congruence as a bundle of
(projective light cones in)
$(2,1)$-planes in $\R^{4,1}$,
cf [\:imdg:, \S6.4.12 and \S6.6.5].
Similarly, a $2$-sphere can be encoded by a $(3,1)$-plane
in $\R^{4,1}$ or, equivalently, by its orthogonal complement:
normalizing yields an identification of hyperspheres with
points in the Lorentz sphere
$$
  S^{3,1} = \{y\in\R^{4,1}\,|\,(y,y)=1\},
$$
where a choice of sign can be used to encode orientation
of a hypersphere.
In particular we obtain representatives in $S^{3,1}$ for
a sphere with centre $m$ and radius $r$ or
a plane with unit normal $n$ through $x$ in $\R^3$
by (cf [\:imdg:, Sect 1.1])
$$
  s = {1\over r}(o+m+{(m,m)-r^2\over 2}\,q)
   \enspace{\rm resp}\enspace
  t = n+(n,x)\,q.
$$

Note that the above definition of Ribaucour pairs of curves
resonates well with the discrete case,
cf [\:boje01:, Sect 4] or [\:imdg:, \S8.3.16],
and \:fig.ribd::
two discrete curves $x,\hat x:\Z\supset I\to\R^3$
are said to form a {\em Ribaucour pair\/} if endpoints of
corresponding edges are concircular,
that is, for adjacent $i,j\in I$,
$$
  c_{ij} :
  = \langle\xi_i,\xi_j,\hat\xi_i,\hat\xi_j\rangle
  \cong \R^{2,1}\subset\R^{4,1}.
\eqno ecircle$$

 \pdfximage width 200pt {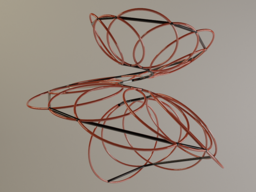}
  \xdef\ribdA{\the\pdflastximage}
 \pdfximage width 200pt {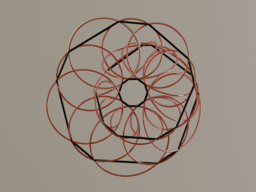}
  \xdef\ribdB{\the\pdflastximage}
 \vglue 0pt\hbox to \hsize{\hfil
  \vtop{\pdfrefximage\ribdA}%
  \hfil
  \vtop{\pdfrefximage\ribdB}%
 \hfil}\vglue 6pt\centerline{%
  {\bf\@fig.ribd:Fig 2@.}
   A Ribaucour pair of discrete curves}\par

Now let $n:I\to S^2$ denote a unit normal field along
one curve $x:I\to\R^3$ of a Ribaucour pair:
this uniquely defines a congruence of hyperspheres
$s(u)$ that contain corresponding points $x(u)$
and $\hat x(u)$ of both curves and have $n(u)$ as a normal;
aligning orientations this construction yields an isometric
isomorphism of normal bundles for the curves of a Ribaucour
pair.
Indeed, employing M\"obius geometric lifts
$$
  t = n + (n,x)\,q \in S^{3,1}
\eqno plift$$
of unit normals, the corresponding (hyper-)spheres can be
encoded as a map into the Lorentz sphere,
$$
  s = t + {1\over r}\,\xi:I\to S^{3,1}
   \enspace{\rm with}\enspace
  {1\over r} =  -{(\hat\xi,t)\over(\hat\xi,\xi)},
\eqno csphere$$
and the induced normal field $\hat n$ of $\hat x$ is
hence given by
$$
  \hat t
  = s - {1\over r}\,\hat\xi
  = t + {1\over r}(\xi-\hat\xi):
  I\to S^{3,1}.
\eqno nisomc$$

\proclaim\@thm.pnf1:Lemma 2.1@.
Suppose $x,\hat x:I\to\R^3$ form a Ribaucour pair of curves;
then the isomorphism {\rm(\:nisomc:)} of normal bundles maps
parallel normal fields $n$ of $x$ to parallel normal fields
$\hat n$ of $\hat x$,
that is, {\rm(\:nisomc:)} intertwines normal connections
of the two curves.

To substantiate this claim first observe that a normal field
$n$ is parallel if and only if its lift (\:plift:)
is,
$$
  n'\parallel x'
   \enspace\Leftrightarrow\enspace
  t'\parallel\xi',
$$
hence if and only if any sphere congruence
$$
  s = t + {1\over r}\,\xi
   \enspace{\rm satisfies}\enspace
  s'\in\Gamma(\langle\xi,\xi'\rangle).
\eqno psphere$$
On the other hand, the enveloped sphere congruence (\:csphere:)
satisfies
$$
  s = t + {1\over r}\,\xi
    = \hat t + {1\over r}\,\hat\xi
  \perp c = \langle\xi,\xi',\hat\xi,\hat\xi'\rangle;
$$
consequently, $s'\perp\xi,\hat\xi$ so that
the claim follows:
$s'\in\Gamma(\langle\xi,\xi'\rangle)$
if and only if
$s'\in\Gamma(\langle\hat\xi,\hat\xi'\rangle)$.

On the other hand, reversing this line of thought, we are able
to reconstruct the enveloped circle congruence of a Ribaucour
pair from a suitable enveloped (hyper-)sphere congruence:
suppose that a sphere congruence $s:I\to S^{3,1}$ touches
two curves $x$ and $\hat x$ and yields parallel normal fields
for both curves,
$$
  s \perp \xi,\xi',\hat\xi,\hat\xi'
   \enspace{\rm and}\enspace
  s'
   \in \Gamma(\langle\xi,\xi'\rangle)
   \cap \Gamma(\langle\hat\xi,\hat\xi'\rangle);
$$
as $s'\perp\xi,\hat\xi$ we conclude that
$$\left.\matrix{
  \xi' \in \Gamma(\langle\xi,s'\rangle) \cr
  \hat\xi' \in \Gamma(\langle\hat\xi,s'\rangle) \cr
  }\right\}\enspace\Rightarrow\enspace
  \langle\xi,\xi',\hat\xi,\hat\xi'\rangle
  = \langle\xi,\hat\xi,s'\rangle.
$$
Thus we have proved:

\proclaim\@thm.pnf2:Lemma 2.2@.
Suppose that a curve $s:I\to S^{3,1}$ of spheres touches
two curves $x,\hat x:I\to\R^3$ and yields parallel normal
fields for these curves;
then $x$ and $\hat x$ form a Ribaucour pair of curves.

As a curve on a surface is a curvature line if and only if
the Gauss map of the surface yields a parallel normal field
along the curve we deduce the following theorem for the channel
surface enveloping one of the curves of spheres discussed
above,
cf [\:pesz17:, Sect 5]:

\proclaim\@thm.rch:Cor 2.3@.
Any two (non-circular) curvature lines of a channel surface
form a Ribaucour pair of curves;
given a Ribaucour pair of curves $x,\hat x:I\to\R^3$ there is
a $1$-parameter family of channel surfaces that contain both
curves as curvature lines.

The first claim follows directly from \:thm.pnf2:,
as the $1$-parameter family of spheres enveloped by a channel
surface yields parallel normal fields along any two of its
(non-circular) curvature lines;
the second claim follows from \:thm.pnf1:,
by using one of the sphere congruences (\:csphere:) given
by a parallel normal field $n$ along $x$ to obtain one of
the sought-after channel surfaces.

Recall that a semi-discrete curvature line net can be thought
of as a sequence of curves, where subsequent curves form
Ribaucour pairs,
cf [\:muwa13:, Def 1.1] or [\:bjmr16:, Sect 3].
As the construction of \:thm.rch: of a channel surface
from a Ribaucour pair of curves involves the choice of
a parallel normal field along one of the two curves,
this construction may be iterated to ``smoothen''
a semi-discrete curvature line net,
see \:fig.chls:.
reminiscent of the smoothing of fully discrete curvature line
nets by using Dupin cyclides,
cf [\:bhv12:]:

\proclaim\@thm.sch:Cor 2.4@.
Any semi-discrete curvature line net can be ``smoothed'' by
a sequence of channel surfaces:
it lies on a $C^1$-surface composed of channel surfaces that
meet at the curves of the semi-discrete net and have the same
tangent planes there.

 \pdfximage width 300pt {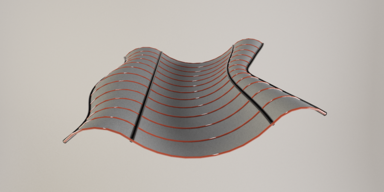}
  \xdef\chls{\the\pdflastximage}
 \vglue 0pt\hbox to \hsize{\hfil
  \vtop{\pdfrefximage\chls}%
 \hfil}\vglue 6pt\centerline{%
  {\bf\@fig.chls:Fig 3@.}
   A semi-discrete curvature line net,
   ``smoothed'' by channel surfaces}\par

This relation between Ribaucour pairs of curves and
enveloped (hyper-)sphere congruences also yields an approach
to reduce the ambient dimension of a curve by means of a
Ribaucour transformation:
given a curve $x:I\to\R^3$ we may use a parallel normal field
$n:I\to S^2$ to construct a sphere congruence touching a given,
fixed sphere and thereby producing a Ribaucour transform of $x$
that is contained in this sphere
---
as the radial vector field of a fixed sphere is parallel along
any curve in that sphere the claim is an immediate consequence
of \:thm.pnf2:.

Using M\"obius geometric techniques again, the enveloped
$1$-parameter family of curves and the sought-after Ribaucour
transform $\hat x$ of $x$ can be determined explicitly, by
algebra alone:
if
$$
  \xi:I\to{\cal L}^4
   \enspace{\rm and}\enspace
  t:I\to S^{3,1}
$$
denote (M\"obius geometric) lifts of the curve $x$ and
a parallel normal field $n$ along $x$, then
$$
  s=t+{1\over r}\,\xi
$$
defines a sphere curve that touches the curve,
$s\perp\xi,\xi'$,
as well as a given fixed sphere $e\in S^{3,1}$
when we let ${1\over r}:={1-(e,t)\over(e,\xi)}$
so that $(s,e)\equiv 1$;
$$
  \hat\xi := s-e:
   I\to{\cal L}^4
\eqno drec$$
then yields (a light cone lift of) the touching point
of $s$ and $e$,
that is,
of the desired Ribaucour transform $\hat x$ of $x$.

\proclaim\@thm.sdr:Cor 2.5@.
Given a (fixed) sphere, any curve $x:I\to\R^3$ that does not
meet the sphere can be Ribaucour transformed into the sphere
by means of a parallel normal field $n$ along $x$;
once a parallel normal field is given, the construction of
the transform is algebraic.

 \pdfximage width 200pt {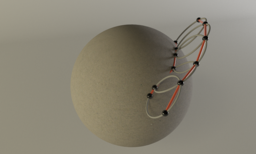}
  \xdef\redA{\the\pdflastximage}
 \pdfximage width 200pt {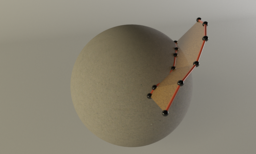}
  \xdef\redB{\the\pdflastximage}
 \vglue 0pt\hbox to \hsize{\hfil
  \vtop{\pdfrefximage\redA}%
  \hfil
  \vtop{\pdfrefximage\redB}%
 \hfil}\vglue 6pt\centerline{%
  {\bf\@fig.redd:Fig 4@.}
   Ambient dimension reduction by means
   or Ribaucour transformation}\par

The first part of \:thm.sdr: works in a completely
analogous way for discrete curves:
given a sphere that the curve does not meet,
a Ribaucour transform of the curve can be constructed
by iteratively constructing second intersection points of
the sphere with the circumcircles of the endpoints of an edge
and an initial point for the edge on the sphere,
cf \:fig.redd::

\proclaim\@thm.ddr:Thm 2.6@.
Given a sphere, any discrete curve $x:\Z\supset I\to\R^3$
that does not meet the sphere can be Ribaucour transformed
onto the sphere;
this construction depends on the choice of one initial point
on the sphere.

At no point of the argument did we use that the ambient
dimension of the curve $x$ be $3$, hence \:thm.sdr:
and \:thm.ddr: hold in any dimension
---
in particular, also for planar curves:
this leads to {\em Ribaucour coordinates\/} for a given curve,
similar to the Ribaucour coordinates used by the classical
geometers for surfaces, cf [\:liro02:],
and generalized for submanifolds in [\:dft07:, Thm 1].
Namely, iterating the dimension reduction of \:thm.sdr:
(or \:thm.ddr:) we obtain a circular arc as a double
Ribaucour transform of a space curve $x:I\to\R^3$.
Suitably aligning iterated transformations,
see \:thm.src: and \:thm.drc:,
the order of transformations is irrelevant and we obtain
Ribaucour coordinates for the given curve:

\proclaim\@thm.rc1:Cor \& Def 2.7@.
Any curve $x:I\to\R^3$ is locally obtained by two subsequent
commuting Ribaucour transforms from a circular arc;\\
the coordinates for the curve obtained in this way
are called {\em Ribaucour coordinates\/} of the curve.

A similar line of thought can be used to interpolate between
curves by Ribaucour transformations.
In the discrete case the construction is straightforward:
once the ambient dimension of the curves has been reduced
to obtain a pair of spherical (or planar) discrete curves,
an ``interpolating'' Ribaucour transform of both curves may
be constructed iteratively, by intersecting circumcircles of
corresponding edges and a common initial point for the edges,
cf \:fig.intd:.

 \pdfximage width 200pt {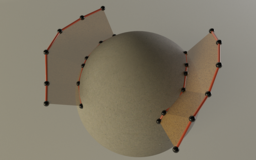}
  \xdef\intA{\the\pdflastximage}
 \pdfximage width 200pt {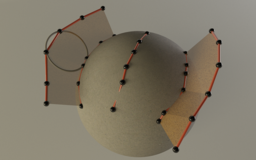}
  \xdef\intB{\the\pdflastximage}
 \pdfximage width 200pt {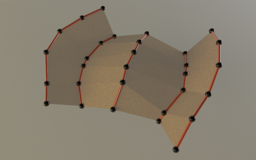}
  \xdef\intC{\the\pdflastximage}
 \pdfximage width 200pt {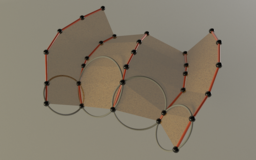}
  \xdef\intD{\the\pdflastximage}
 \vglue 0pt\hbox to \hsize{\hfil
  \vtop{\pdfrefximage\intA}%
  \hfil
  \vtop{\pdfrefximage\intB}%
 \hfil}\vglue 6pt\hbox to \hsize{\hfil
  \vtop{\pdfrefximage\intC}%
  \hfil
  \vtop{\pdfrefximage\intD}%
 \hfil}\vglue 6pt\centerline{%
  {\bf\@fig.intd:Fig 5@.}
   Interpolation by Ribaucour transformations}\par

In the smooth case the argument is similar:
after a dimension reduction a common Ribaucour transform
of two planar or spherical curves needs to be constructed.
This is slightly more involved: a simple existence argument
relies on [\:buje06:, Thm 3.4],
where we use methods from Lie sphere geometry.
Thus we add another (timelike) dimension to the coordinate
space of the M\"obius geometric model,
$$
  \R^{3,2}
  = \langle p\rangle\oplus_\perp\R^{3,1}
  = \langle o,p,q\rangle\oplus_\perp\R^2,
$$
cf [\:bl29:] or [\:buje06:];
we then get a {\em Legendre lift\/} of a curve $x:I\to\R^2$
with unit normal field $n:I\to S^1$ by
$$
  \lambda := \langle\xi,\nu\rangle
   \enspace{\rm with}\enspace
  \nu := p+t = p+n+(n,x)\,q,
$$
that is,
a curve of null $2$-planes in $\R^{3,2}$, resp lines
in the Lie quadric $\P({\cal L}^4)\subset\P(\R^{3,2})$,
that satisfies the contact condition,
$$
  \sigma'\perp\lambda
   \enspace{\rm for}\enspace
  \sigma\in\Gamma(\lambda).
$$
Now suppose that two curves $x_i:I\to\R^2$ ($i=0,1$) with
unit normal fields $n_i$ are given and
assume, for regularity, that the two curves do not touch
a common circle at any two corresponding points $x_i(u)$.
Then
$$
  b := \lambda_0\oplus\lambda_1:I\to G_{2,2}(\R^{3,2})
\eqno demoulin$$
defines a $(2,2)$-bundle as a map into the Grassmannian of
$(2,2)$-planes in $\R^{3,2}$.
As this is a bundle over a $1$-dimensional base it is clearly
a flat bundle,
thus [\:buje06:, Thm 3.4] assures that any null $2$-plane
subbundle $\hat\lambda$ of $b$ that intersects $\lambda_i$
non-trivially defines a common Ribaucour transform
$\hat x:I\to\R^2$ of the two initial curves $x_i$.

Thus we obtain the following theorem that holds for both,
smooth as well as discrete curves:

\proclaim\@thm.ri2:Thm 2.8@.
Any two planar curves $x_0,x_1:I\to\R^2$ admit,
under mild regularity assumptions,
a common Ribaucour transform $\hat x:I\to\R^2$;
in the smooth and discrete cases there is a $1$- resp
$2$-parameter family of such common Ribaucour transforms.

As a consequence any two (discrete or smooth) space curves
span a discrete resp semi-discrete curvature line net:
we can solve a mixed boundary and initial value problem
for (semi-)discrete curvature line nets.
More precisely:

\proclaim\@thm.ri3:Cor 2.9@.
Any two discrete or smooth space curves $x_0,x_1:I\to\R^3$
form two boundary curves of a discrete resp semi-discrete
curvature line net;
generically, three interpolating Ribaucour transforms are
required.

\section 3. Submanifolds

After putting forward the key ideas of this work in the
previous section we are now prepared to tackle the general
case of submanifolds in $\R^n$ or $S^n$.
In the case of space curves, we based our presentation on
two lemmas to then easily deduce the results we were
interested in
---
however, while \:thm.pnf1: is insensitive to dimensions,
and a useful reformulation will be readily available,
\:thm.pnf2: crucially depends on the ``correct'' dimensions,
and we will therefore need to use different arguments to prove
the results that depend on it.

In fact, the notion of a Ribaucour pair resp transform
is more intricate in higher dimensions:
in addition to the fact that two submanifolds envelop
a congruence of spheres of the right dimension,
a generalization of the classical demand that
``curvature lines correspond''
on the two envelopes needs to be implemented,
cf [\:imdg:, Def 8.2.2] and [\:buje06:, Sect 4]:

\proclaim\@def.ribs:Def 3.1@.
Two immersed submanifolds
$x,\hat x:\Sigma^m\to S^n=\P({\cal L}^{n+1})$
in the conformal $n$-sphere,
where ${\cal L}^{n+1}=\{y\in\R^{n+1,1}\,|\,(y,y)=0\}$,
form a {\em Ribaucour pair\/} if
{\parindent=2em
\item{\rm(i)} they envelop a congruence of $m$-spheres;
\item{\rm(ii)} the $(1,1)$-subbundle $x\oplus\hat x$
 in $\R^{n+1,1}$ is flat.}

Note that \:def.ribs:(ii) yields a weak version of the
condition that the curvature directions of two hypersurfaces
correspond.
Also note the slight change of notation:
working in the conformal $n$-sphere from the start,
$x$ and $\hat x$ no longer denote immersions into $\R^n$
but into $S^n=\P({\cal L}^{n+1})$.

\proclaim\@thm.ribt:Lemma 3.2@.
Two pointwise distinct immersions
$x,\hat x:\Sigma^m\to S^n=\P({\cal L}^{n+1})$
{\parindent=2em
\item{\rm(i)} envelop an $m$-sphere congruence
 if and only if $\fdeb x\oplus\hat x=x\oplus\fdeb{\hat x}$,
 where we denote
 $$
   \fdeb x|_u:=x(u)\oplus d_u\xi(T_u\Sigma^m)
    \enspace{\sl for}\enspace
   x=\langle\xi\rangle
     \enspace{\sl and}\enspace
   u\in\Sigma^m;
 \eqno derb$$
\item{\rm(ii)} form a Ribaucour pair
 if and only if $d\hat\xi\in\Omega^1(\fdeb x)$
 for a suitable lift $\hat\xi\in\Gamma(\hat x)$ of $\hat x$.}

\:thm.ribt:(i) says that the $m$-spheres touching $x$ and
containing the points of $\hat x$ coincide with those that
touch $\hat x$ and contain the corresponding points of $x$;
\:thm.ribt:(ii) implies \:thm.ribt:(i) and says that
$d\hat\xi\perp x$,
hence that there is a parallel section $\hat\xi$
of $x\oplus\hat x$.

At some points it will be useful to consider Legendre
lifts into Lie sphere geometry,
cf [\:buje06:],
however we will mostly stay in M\"obius geometry
as we aim to keep the notions of points and circles,
in particular with a view to the discrete setting,
cf [\:boje01:] or [\:imdg:, \S8.3.16]:

\proclaim\@def.ribd:Def 3.3@.
Two discrete circular nets
$x,\hat x:\Z^m\supset\Sigma^m\to S^n=\P({\cal L}^{n+1})$
in the conformal $n$-sphere
form a {\em Ribaucour pair\/} if
corresponding edges are {\em circular\/},
that is, have concircular endpoints.

If two immersed submanifolds $x,\hat x:\Sigma^m\to S^n$
envelop an $m$-sphere congruence then we obtain a natural
identification
$$
  \wlnb x
  \cong (\fdeb x\oplus\hat x)^\perp
  \cong (x\oplus\fdeb{\hat x})^\perp
  \cong \wlnb{\hat x}
$$
of their (weightless) normal bundles,
cf [\:buje06:, (4.4)] or [\:buca10:, Sect 8.2]:
geometrically, this is the identification of \:thm.pnf1:,
where a normal of $x$ is mapped to a normal of $\hat x$
by means of a hypersphere that contains the $m$-sphere
of the enveloped sphere congruence,
$$
  \hat t - {(\xi,\hat t)\over(\xi,\hat\xi)}\,\hat\xi
  = s
  = t - {(\hat\xi,t)\over(\hat\xi,\xi)}\,\xi
   \enspace{\rm hence}\enspace\cases{
  \hat t + \hat x = s + \hat x & and \cr
  t + x = s + x. & \cr}
\eqno nisomm$$
Hence we obtain a higher dimensional version of \:thm.pnf1::

\proclaim\@thm.pnfm:Lemma 3.4@.
If two immersed submanifolds $x,\hat x:\Sigma^m\to S^n$
envelop an $m$-sphere congruence
then there is a natural, connection preserving isometric
isomorphism {\rm(\:nisomm:)} between their (weightless)
normal bundles $\wlnb x$ and $\wlnb{\hat x}$.

A higher dimensional version of the second claim of \:thm.rch:
will be obtained as an immediate consequence,
as in the case of space curves:
any parallel normal field of one envelope of an $m$-sphere
congruence yields one for the other via a ``parallel''
enveloped hypersphere congruence $s$,
hence both envelopes form {\em extended curvature leaves\/}
of the envelope of the $m$-parameter family of hyperspheres,
in the sense that their tangent spaces are invariant subspaces
for the shape operator of the hypersurface,
since
$$
  0 = \nabla^\perp(t+x)
   \enspace\Leftrightarrow\enspace
  d(t+x) \in \Omega^1(\fdeb x/x),
$$
and similarly for the normal field $\hat t+\hat x$ that
corresponds to $t+x$ via (\:nisomm:).
Conversely, the envelope of a (spacelike) $m$-parameter family
of hyperspheres $s:\Sigma^m\to S^{n,1}$ is foliated by
spherical curvature leaves
$$
  c(u) := (s(u)\oplus d_us(T_u\Sigma^m))^\perp
   \cong \R^{n-m,1} \subset \R^{n+1,1}
$$
of dimension $\dim c(u)-2=n-m-1$;
as long as $m<n-1$ any two extended curvature leaves
$x$ and $\hat x$ complementary to the spherical curvature
leaves $c$
are tangent to the $m$-spheres
$$
  (x\oplus\hat x\oplus c^\perp)(u)\cap s^\perp(u)
  = x(u)\oplus\hat x(u)\oplus d_us(T_u\Sigma^m)
   \cong \R^{m+1,1} \subset \R^{n+1,1}.
$$
Note that, in the case $m=n-1$,
the spherical ``curvature leaves'' degenerate to point pairs
and we obtain the familiar figure of a hypersphere congruence
with its envelope consisting of two sheets.

Summarizing, we obtain the following higher dimensional
version of \:thm.rch::

\proclaim\@thm.rchm:Cor 3.5@.
If $x,\hat x:\Sigma^m\to S^n$ envelop an $m$-sphere congruence,
and if $x$ has a parallel normal field,
then $x$ and $\hat x$ yield extended curvature leaves
of a hypersurface that is obtained as the envelope
of an $m$-parameter family of hyperspheres.\\
Conversely, given the envelope of an $m$-parameter family
$s:\Sigma^m\to S^{n,1}$ of hyperspheres
in $S^n$, where $m<n-1$, any two extended curvature leaves
complementary to its spherical curvature leaves envelop
a congruence of $m$-spheres.

This generalization of \:thm.rch: shows that the construction
of a channel surface from a Ribaucour pair of curves only
hinges on the enveloped circle congruence and the
existence of parallel normal fields,
not on the fact that they form a ``Ribaucour pair''
in the more general sense:
if we demand higher dimensional submanifolds $x$ and $\hat x$
to form a Ribaucour pair in the sense of \:def.ribs: then
more fine structure of $x$ and $\hat x$ as extended
curvature leaves of a hypersurface can be derived.

Note that the first construction of \:thm.rchm:
can be iterated when a sequence of submanifolds that
envelop $m$-sphere congruences in a suitable manner
is given and if one of the submanifolds admits a
parallel normal field:
this yields a higher dimensional version of \:thm.sch:.

To get a higher dimensional version of the dimension reduction
procedure \:thm.sdr: we follow the same line of thought as
before:
given a submanifold $x:\Sigma^m\to S^n$
and a fixed hypersphere $e\in S^{n,1}$,
we use a unit normal field $t+x\in\Gamma(\wlnb x)$
and a lift $\xi\in\Gamma(x)$ of $x$
to construct a congruence of hyperspheres
$$
  s = t + {1-(e,t)\over(e,\xi)}\,\xi
$$
that simultaneously touches $x$ and $e$, that is,
$s\perp\xi,d\xi$ and $(s,e)\equiv 1$.
The touching points
$$
  \hat x = \langle\hat\xi\rangle
   \enspace{\rm with}\enspace
  \hat\xi := s-e
$$
of $s$ and $e$ then form a Ribaucour transform of $x$,
by \:thm.ribt:,
as soon as $t+x$ is a parallel normal field of $x$,
since
$$
  d\hat\xi = ds \equiv dt \bmod \fdeb x;
$$
assuming that $\hat x$ is complementary to $x$,
that is, that $x$ does not intersect the hypersphere $e$,
we learn that this condition is also necessary:

\proclaim\@thm.sdrm:Cor 3.6@.
A submanifold $x:\Sigma^m\to S^n$ admits a Ribaucour transform
to a hypersphere $e\in S^{n,1}$ that it does not meet,
$x\not\perp e$,
if and only if
$x$ admits a parallel normal field $t+x\in\Gamma(\wlnb x)$;
in this case a spherical Ribaucour transform is given by
$$
  \hat x = \langle s-e\rangle
   \enspace{\rm with}\enspace
  s := t + {1-(e,t)\over(e,\xi)}\,\xi
   \enspace{\rm and}\enspace
  \xi \in \Gamma(x).
\eqno drem$$

Note that, in this dimension reduction construction,
$x$ and $\hat x$ form a Ribaucour pair as soon as
they envelop an $m$-sphere congruence:
since $(ds,\xi)=-(s,d\xi)=0$ we conclude that
$$
  ds \in \Omega^1(\fdeb x\oplus\hat x)
   \enspace\Leftrightarrow\enspace
  ds \in \Omega^1(\fdeb x).
$$

This observation is confirmed by the discrete case.
Given a circular net
$
  x:\Z^m\supset\Sigma^m\to S^n,
$
and a hypersphere $e\in S^{n,1}$ that the net does not meet
we construct a spherical net iteratively,
by evolving along edges using the intersection of the fixed
sphere with the circumcircle of an edge and an initial point
for the edge, as in the $1$-dimensional case \:thm.ddr:.
This construction is consistent, and yields a circular net,
by the M\"obius geometric formulation of Miguel's theorem:
{\em given seven vertices of a cube with circular faces,
 the eighth vertex can be constructed uniquely\/}.
Namely, considering a face of the original net $x$ and
an initial point of the corresponding face of $\hat x$
the construction of the face of $\hat x$ is completed
on the $2$-sphere containing the five starting points,
cf [\:boje01:, Sect 4].
A standard dimension count then shows that corresponding
$m$-cells of the circular nets $x$ and $\hat x$ lie on
$m$-spheres:

\proclaim\@thm.ddrm:Cor 3.7@.
Any circular net $x:\Z^m\supset\Sigma^m\to S^n$ admits
a Ribaucour transform to a fixed hypersphere $e\in S^{n,1}$
that it does not meet, $x\not\perp e$;
the two nets {\em envelop\/} a discrete $m$-sphere congruence,
in the sense that corresponding $m$-cells lie on $m$-spheres.

As in the case of curves in \:thm.rc1:,
the dimension reductions of \:thm.sdrm: resp \:thm.ddrm:
can be used to introduce Ribaucour coordinates
for a submanifold with flat normal bundle
---
where ``enough'' spherical Ribaucour transforms exist,
cf [\:dft07:, Thm 1].

Thus suppose that $t_i+x\in\Gamma(\wlnb x)$, $i=1,\dots,n-m$,
are parallel unit normal fields of an immersed submanifold
$x:\Sigma^m\to S^n$ and
let $e_1,\dots,e_{n-m}\in S^{n,1}$ denote an orthonormal
system that defines an $m$-sphere in the conformal $n$-sphere
via
$$
  e = \langle e_i\,|\,i=1,\dots,m\rangle^\perp
   \cong \R^{m+1,1} \subset \R^{n+1,1}.
$$
Now we let $\xi\in\Gamma(x)$ and iterate (\:drem:) to obtain
an iterated Ribaucour transform of $x$ into the $m$-sphere $e$:
with a first Ribaucour transformation of $x$ into the $i$-th
hypersphere $e_i$ and its parallel normal field obtained by
(\:nisomm:) from the $j$-th parallel normal field of $x$,
$$
  \xi_i = (t_i-e_i) + {1-(e_i,t_i)\over(e_i,\xi)}\,\xi
   \enspace{\rm and}\enspace
  t_{ij} - {(\xi,t_{ij})\over(\xi,\xi_i)}\,\xi_i
  = t_j - {(\xi_i,t_j)\over(\xi,\xi_i)}\,\xi,
$$
a lengthy but straightforward computation reveals that
$$
  \xi_{ij}
  = (t_j-e_j)
  + {(1-(e_j,t_j))(e_i,\xi)+(e_i,t_j)(e_j,\xi)\over
     (1-(e_i,t_i))(e_j,\xi)+(e_j,t_i)(e_i,\xi)}\,(t_i-e_i)
  + {(1-(e_i,t_i))(1-(e_j,t_j))-(e_i,t_j)(e_j,t_i)\over
     (1-(e_i,t_i))(e_j,\xi)+(e_j,t_i)(e_i,\xi)}\,\xi;
\eqno permij$$
we observe that, up to scaling, $\xi_{ij}$ is symmetric in
$i$ and $j$, thus confirming Bianchi's permutability theorem
for the particular Ribaucour transformations we use,
cf [\:buje06:, Sect 3] or [\:dft07:, Cor 16].
Once the circularity claim of Bianchi's permutability theorem
is established for $x$, $x_i$, $x_j$ and $x_{ij}$,
$$
  x_{ij} = \langle\xi_{ij}\rangle
   \enspace{\rm with}\enspace
  \xi_{ij} \in\Gamma(x\oplus x_i\oplus x_j),
$$
we obtain an alternative approach to determining $x_{ij}$:
the conditions $\xi_{ij}\perp e_i,e_j$ for $x_{ij}$ lead to
a homogeneous system of linear equations for the coefficients
in (\:permij:),
$$
  \xi_{ij} = a_0\xi + a_i(t_i-e_i) + a_j(t_j-e_j)
   \enspace{\rm with}\enspace\cases{
   0 
    = (e_i,\xi)\,a_0 - (1-(e_i,t_i))\,a_i + (e_i,t_j)\,a_j,\cr
   0 
    = (e_j,\xi)\,a_0 + (e_j,t_i)\,a_i - (1-(e_j,t_j))\,a_j;\cr}
\eqno bquad$$
Cramer's rule then recovers the coefficients of (\:permij:)
up to scaling --- as long as the equations are independent,
which can be achieved by ensuring that
$x$ does not hit the hyperspheres $e_i^\perp$ and $e_j^\perp$,
and $x_i$ and $x_j$ do not hit the target $(n-2)$-sphere
$\langle e_i,e_j\rangle^\perp$ of $x_{ij}$:
we require that the dimension reduction of \:thm.sdrm: works
in every step of the iteration.

Using (\:bquad:) it is now straightforward to formulate the
general case:
as every permutation of a finite set is a composition of
transpositions our dimension reduction is independent
of order,
and circularity of Bianchi quadrilaterals ensures that
(\:bquad:) generalizes.
More precisely,

\proclaim\@thm.src:Thm \& Def 3.8@.
Suppose that $t_i+x\in\Gamma(\wlnb x)$, $i=1,\dots,n-m$,
are parallel unit normal fields of an immersed submanifold
$x:\Sigma^m\to S^n$ and
let
$
  e = \langle e_i\,|\,i=1,\dots,m\rangle^\perp
$
denote an $m$-sphere, given in terms of orthogonally
intersecting hyperspheres $e_1,\dots,e_{n-m}\in S^{n,1}$;
then
$$
  \hat\xi := a_0\xi + \sum_{i=1}^{n-m}a_i\,(t_i-e_i) \neq 0
   \enspace{\rm with}\enspace
  (\hat\xi,e_i) = 0
   \enspace{\rm for}\enspace
  i=1,\dots,n-m
$$
defines an $(n-m)$-fold Ribaucour transform
$
  \hat x = \langle\hat\xi\rangle:\Sigma^m\to S^n
$
of $x$ into the $m$-sphere $e$.\\
We say that $\hat x$ defines {\em Ribaucour coordinates\/}
for the submanifold $x$.

Note that the choices of (parallel resp constant) normal fields
for $x$ and the $m$-sphere $e$ establishes an isometric and
connection preserving isomorphism between normal bundles:
this isomorphism is the source of our vectorial Ribaucour
transformation that establishes the Ribaucour coordinates,
cf [\:dft07:, Thm 1].

Based on \:thm.ddrm: Ribaucour coordinates for discrete
principal nets can be introduced in a similar way:
in order to ensure that permutability holds we need to choose
initial points for the spherical Ribaucour transforms suitably
---
Miguel's theorem then guarantees that the double Ribaucour
reduction $x_{ij}$ is independent of order and corresponding
points of the circular nets $x$, $x_i$, $x_j$ and $x_{ij}$ are
concircular,
see [\:imdg:, Sect 8.5].
Thus in the case of an $m$-dimensional net in $S^n$ initial
points need to be suitably chosen on an $(n-m-1)$-sphere
containing the initial point of $x$ and the corresponding
point of the target net defining the Ribaucour coordinates.

\proclaim\@thm.drc:Thm \& Def 3.9@.
Suppose that $x:\Z^m\supset\Sigma^m\to S^n$ is a circular net
and fix an $m$-sphere
$
  e = \langle e_i\,|\,i=1,\dots,m\rangle^\perp
$
in terms of orthogonally intersecting hyperspheres
$e_1,\dots,e_{n-m}\in S^{n,1}$;
further fix $i_0\in\Sigma^m$ and
choose an $(n-m)$-cube
$$
  \{0,1\}^{n-m}\ni\varepsilon \mapsto y_\varepsilon\in S^n
   \enspace{\rm with}\enspace\cases{
   y_\varepsilon = x(i_0) & for $\varepsilon=0$ and \cr
   y_\varepsilon \in \langle\varepsilon_1e_1,\dots,
    \varepsilon_{n-m}e_{n-m}\rangle^\perp &
    for $\varepsilon\neq 0$ \cr}
$$
and with circular faces.
Then the $(n-m)$-fold iterated Ribaucour reduction
with initial points $y_\varepsilon$
$$
  \hat x:\Sigma^m\to e =
   \langle e_1,\dots,e_{n-m}\rangle^\perp\subset S^n,
$$
and is independent of the order of transformations.\\
We say that $\hat x$ defines {\em Ribaucour coordinates\/}
for the discrete circular net $x$.

While the dimension reduction via Ribaucour transformations
works perfectly in higher dimensions,
the interpolation by Ribaucour transformations of \:thm.ri2:
hinges on the flatness of the Demoulin vector bundle
(\:demoulin:)
---
which is no longer trivially satisfied in higher dimensions.
Indeed, in [\:buje06:, Sect 6] we provide an example of two
surfaces that have one common (totally umbilic) Ribaucour
transform but not more, so that permutability fails:
in particular,
the curvature lines on the common Ribaucour transform
(obtained from the respective Ribaucour transformation)
do not line up.

This situation is typical:
by introducing Ribaucour coordinates of \:thm.src:
for two $m$-dimensional submanifolds with flat normal bundle
in $S^n$,
the two submanifolds may be (locally) obtained by a sequence
of $2(n-m)$ Ribaucour transformations from each other,
however, their curvature directions will in general
not line up in the desired way.

In the discrete case,
where curvature line coordinates are inherent,
we learn from these observations that interpolation
by Ribaucour transformations fails in general:
a simple discretization of the Dupin cyclides used
in [\:buje06:, Sect 6], by sampling curvature lines,
provides counter-examples.

\section References

\message{References}
\bgroup\frenchspacing\parindent=2em

\refitem bl29
 W Blaschke:
 {\it Vorlesungen \"uber Differentialgeometrie III\/};
 Springer Grundlehren XXIX, Berlin (1929)

\refitem bkpww11
 P Bo, M Kilian, H Pottmann, J Wallner, W Wang:
 {\it Circular arc structures\/};
 ACM SIGGRAPH 101 (2011)

\refitem boje01
 A Bobenko, U Hertrich-Jeromin:
 {\it Orthogonal nets and Clifford algebras\/};
 T\^ohoku Math Publ 20, 7--22 (2001)

\refitem bhv12
 A Bobenko, E Hunen-Venedy:
 {\it Curvature line parametrized surfaces and orthogonal
  coordinate systems:
  Discretization with Dupin cyclides\/};
 Geom Dedicata 159, 207--237 (2012)

\refitem buje06
 F Burstall, U Hertrich-Jeromin:
 {\it The Ribaucour transformation in Lie sphere geometry\/};
 Differ Geom Appl 24, 503--520 (2006)

\refitem buca10
 F Burstall, D Calderbank:
 {\it Conformal submanifold geometry I-III\/};
 EPrint arXiv:1006.5700 (2010)

\refitem bjmr16
 F Burstall, U Hertrich-Jeromin, C M\"uller, W Rossman:
 {\it Semi-discrete isothermic surfaces\/};
 Geom Dedicata 183, 43--58 (2016)

\refitem cote04
 A Corro, K Tenenblat:
 {\it Ribaucour transformations revisited\/};
 Commun Anal Geom 12, 1055--1082 (2004)

\refitem dft07
 M Dajczer, L Florit, R Tojeiro:
 {\it The vectorial Ribaucour transformation for submanifolds
  and its applications\/};
 Trans Amer Math Soc 359, 4977--4997 (2007)

\refitem imdg
 U Hertrich-Jeromin:
 {\it Introduction to M\"obius differential geometry\/};
 London Math Soc Lect Note Series 300,
  Cambridge Univ Press, Cambridge (2003)

\refitem rosz18
 U Hertrich-Jeromin, W Rossman, G Szewieczek:
 {\it Discrete channel surfaces\/};
 in preparation

\refitem liro02
 L de Lima, P Roitman:
 {\it Constant mean curvature one surfaces in hyperbolic
  $3$-space using the Bianchi-Cal\`o method\/};
 Ann Braz Acad Sci 74, 19--24 (2002)

\refitem muwa13
 C M\"uller, J Wallner:
 {\it Semi-discrete isothermic surfaces\/};
 Res Math 63, 1395--1407 (2013)

\refitem pesz17
 M Pember, G Szewieczek:
 {\it Channel surfaces in Lie sphere geometry\/};
 EPrint arXiv:1709.02224 (2017)

\egroup
\vskip3em\vfill
\bgroup\fn[cmr7]\baselineskip=8pt
\def\addwd{\hsize=.36\hsize}
\def\fran{\vtop{\addwd
 F Burstall\\
 Department of Mathematical Sciences\\
 University of Bath\\
 Bath, BA2 7AY (UK)\\
 Email: feb@maths.bath.ac.uk
 }}
\def\udomaria{\vtop{\addwd
 U Hertrich-Jeromin and M Lara Miro\\
 Vienna University of Technology\\
 Wiedner Hauptstra\ss{}e 8--10/104\\
 A-1040 Vienna (Austria)\\
 Email: udo.hertrich-jeromin@tuwien.ac.at\\
 \phantom{Email: }mmiro@geometrie.tuwien.ac.at
 }}
\hbox to \hsize{\hfil \fran \hfil \udomaria \hfil}
\egroup
\end{document}                            